
\magnification= 1440
\vsize=24.7true cm
\hsize=17true cm

\
\bigskip
{\hfill  December, 1994
$\qquad$ }

{\hfill  ITP--94--54E$\qquad$  }

\bigskip
\bigskip
\centerline{\bf COVARIANT DIFFERENTIAL AND INTEGRAL CALCULI }
\smallskip
\centerline{\bf FOR LATTICE (l,q)-DEFORMED FIELDS }
\bigskip
\bigskip
\centerline{Anatolij I.~Bugrij \plainfootnote{$^1$}{E-mail address:
 abugrij\@gluk.apc.org }}
\medskip \centerline{\it Bogolyubov
 Institute for Theoretical Physics,}
 \centerline{\it 252142 Kiev-142,
Ukraine}
\bigskip
\centerline{Vladimir Rubtsov
\plainfootnote{$^2$}{E-mail address: volodya\@orphee.polytechnique.fr}
\plainfootnote{$^3$}{Permanent address: ITEP, B.Cheremushkinskaya 25,
117259, Moscow, Russia } }
\medskip
\centerline{\it  Centre de Mathematiques, Ecole
Polytechnique,} \centerline{\it F-91128 Palaiseau,France}
\bigskip \centerline{Vitalij N.~Shadura{$^1$}}
\medskip
\centerline{\it Bogolyubov Institute for
 Theoretical Physics,}
 \centerline{\it 252142 Kiev-142, Ukraine}
\vskip 1.5true cm
\centerline{\bf Abstract}
\medskip
\noindent
Using the Hecke $\hat R$-matrix, we give a definition of the lattice
$(l,q)$-deformed $n$-component boson and Grassmann fields.
Here $l$ is a deformation parameter for the commutation
relations of "values" of these fields in two arbitrary lattice
sites and $q$ is a deformation parameter for $n$-component
$q$-boson or  $q$-Grassmann variable.
In framework of the
Wess-Zumino approach to the noncommutative differential calculus the
commutation relations between differentials and derivatives of these
fields are determined. The $SL_q(n,C)$-invariant generalization of the
Berezin integration for the lattice $n$-component $(l,q)$-\-Grassmann
field is suggested.  We show that the Gaussian functional integral for
 this field is expressed through the $(l,q)$-deformed counterpart of
the Pfaffian.
\vfill\eject



\define\a{\alpha}
\redefine\b{\beta}
\define\e{\varepsilon}
\define\A{{\hat A}}

\define\R{{\hat R}}
\define\CR{{\hat {\Cal R}}}

\redefine\L{{\hat L}}
\redefine\D{{\hat D}}
\redefine\Q{{\hat Q}}
\define\F{{\hat F}}

\redefine\C{{\hat C}}
\redefine\B{{\hat B}}
\define\CG{\hat {\Cal G}}

\define\V{{\Cal V}^{(n)}}

\redefine\P{{\Cal \Psi}^{(n)}}
\redefine\CP{{\hat P}}

\define\si#1{\sigma_{#1}}
\define\n{\nu}
\define\m{\mu}
\redefine\pmm#1{\p^{\mu}_{#1}}
\define\pn#1{\p^{\nu}_{#1}}
\define\w{\wedge}
\define\K{{\Cal \Xi}^{(n)}}

\define\r{\rho}
\define\ca{{\hat {\tenit A}}}
\define\cs{{\hat {\tenit S}}}

\define\caa{\,{\Cal A}\,}


\define\sa{{SL_{\,q}(n,C)}}
\define\scl{{SL_{\,l}(N,C)}}

\define\caq{{{ \hat {\tenit A}}_{\,q}}}
\define\csq{{{ \hat {\tenit S}}_{\,q}}}

\define\s{{SL_{\,q}(n,C)}}
\define\sq{{SL_{\,q}(n)}}
\define\J{{\hat J}}

\define\CF{{\Cal F}}
\define\CJ{{\hat {\Cal J}}}

\define\CI{{\hat {\Cal I}}}

\define\I{{\hat I}}

\define\kt{{\tilde {\xi}}}
\define\vtt{{\tilde v}}
\define\k{\xi}
\define\g{\gamma}
\redefine\d{\delta}
\define\ty{\otimes}
\redefine\tt{\hat \otimes}
\define\my{{\dot  \otimes}}
\define\ka{\k^{\,\alpha}}
\define\kb{\k^{\,\beta}}

\define\kro{\k^{\r}}
\define\kg{\k^{\g}}

\define\rka#1{\k^{\,\alpha}_{#1}}
\define\rkb#1{\k^{\,\beta}_{#1}}
\define\va{v^{\,\alpha}}

\define\pa#1{\psi^{\a}_{#1}}
\define\pb#1{\psi^{\b}_{#1}}
\define\pg#1{\psi^{\g}_{#1}}
\define\pro#1{\psi^{\r}_{#1}}
\define\pd#1{\psi^{\d}_{#1}}
\define\ro#1{\psi^{#1}_{r_1}}
\define\rr#1{\psi^{#1}_{r_2}}
\define\p{\psi}
\define\tp{\tilde \psi}

\define\fa#1{\varphi^{\a}_{#1}}
\define\fb#1{\varphi^{\b}_{#1}}
\define\fg#1{\varphi^{\g}_{#1}}
\define\fro#1{\varphi^{\r}_{#1}}
\define\fd#1{\varphi^{\d}_{#1}}
\define\vf{\varphi}

\define\kr{\k^{\,\rho}}

\define\mr{\,{\R^{\,\a\b}_{\,\g\r}}\,}

\define\jm#1{\k^{\mu}_{#1}}
\define\jn#1{\k^{\nu}_{#1}}

\define\ja#1{\k^{\a}_{#1}}
\define\jd#1{\k^{\d}_{#1}}
\define\jb#1{\k^{\b}_{#1}}
\define\jg#1{\k^{\g}_{#1}}
\define\jro#1{\k^{\r}_{#1}}
\define\cD{{\Cal D}}

\define\iro#1{\partial ^{\,#1} _{\r}}
\define\ia#1{\partial ^{\,#1} _{\alpha}}
\define\ib#1{\partial  ^{\,#1}_{\beta}}
\define\ig#1{\partial  ^{\,#1}_{\gamma}}
\define\id#1{\partial  ^{\,#1}_{\delta}}

\define\baq{{{\hat {\Cal A}}^{(l,q)}}}
\define\bsq{{{\hat {\Cal S}}^{\,(l,q)}}}

\TagsOnRight

\centerline{\bf A.~Bugrij, V.~Rubtsov, V.~Shadura}
\bigskip

\centerline{\bf Covariant Differential and Integral Calculi}
\centerline{\bf for Lattice
(l,q)-Deformed Fields}
\bigskip

\centerline{\bf 1. Introduction}
\bigskip

During the two last decades the interplay between  approaches
to investigation of  lattice field theories (as a method of
regularization of quantum field theory) and  critical phenomena in
the lattice spin systems influenced considerably on
understanding of the non-perturbative phenomena in the  quantum field
theory. On the other hand the possibility of the partition function
representation of the lattice spin systems through partition function
of the lattice field theories allows to
describe  the critical behaviour of these systems by means of a
quantum field theory.  So, for example, in the critical point
  partition function of two-dimensional Ising model can be
 represented as   partition function of the lattice massless majorana
fermion field theory [1-3].
In the continuum limit the masslessness of
the majorana fermion is a consequence of conformal invariance of the
 critical fluctuations.  As it was shown by Belavin, Polyakov and
 Zamolodchikov in [4], this invariance
  allows  to assign to each universality class of the critical
   behaviour of exactly solvable lattice systems
the corresponding two-dimensional quantum conformal field theory. This
 result is connected with the fact that, in critical points, such
 systems, besides the invariance with respect to the global conformal
 transformations ($SL(2,C)$ group) possess larger symmetry, that is,
  the space of eigenvectors of the transfer matrix has an
 ivariance with respect to the infinite-dimensional Virasoro algebra.

 However, as soon as the correlation
 radius of critical fluctuations is finite nearby the critical point,
 the conformal invariance is broken.  Nevertheless, as it was shown
 by Zamolodchikov in [5], an exactly solvable two-dimentional quantum
 field theory can still possess infinitely many integrals of motion,
while being nearby the fixed point (this point corresponds to a quantum
 conformal field theory).
The infinite number  integrals of   motions was  found [6]
near the critical point in the eight-vertex model  in the point where this
model corresponds to the two noninteracting Ising sub-lattices.
Recently Jimbo and Miwa with collaborators showed in the series of papers
[7,8] that the eigenvectors of the corner transfer-matrix
for the six-vertex model form the irreducible representation space
of  $U_q(\hat {sl_2})$ at $-1<q<0$. Moreover, in [9,10]   was shown,
that the eigenvectors space of the
corner transfer-matrix for the RSOS  model and, in particular, for
the two-dimensional Ising model [10] can also possess this symmetry.

These results allow one to assume that the partition function of the
 two-dimensional exactly solvable lattice model can be expressed
through the partition function of some $SL_q(2,C)$-invariant lattice
$q$-deformed field theory, Fock space of which being constructed using
the representation theory of $U_q(\hat {sl_2})$.

The attempt to realize this idea was suggested in [3]
where the  representation of the partition function of the
two-dimensional Ising model in the form of the $SL_q(2,R)$-invariant
functional integral over the lattice real  $q$-fermion field $(q=-1)$ was
found.
The following form of  commutation relations for
the lattice $q$-fermion
field was proposed in [3]:
  $$
\k^{\a}_{i}\,\k^{\b}_{j}=-q\,{\R^{\,\a\b}_{\,\g\r}}\,
{\CP^{\,k\,l}_{\,i\,j}}\,\k^{\g}_{k}\,\k^{\r}_{l}\,\tag1.1
$$
where
$\CP^{\,k\,l}_{\,i\,j}=\d^{k}_{i}\,\d^{l}_{j}+\d^{k}_{j}\,\d^{l}_{i}$
is the permutation matrix
and the values of a latin
indices  goes over all values of variable $r$ which
 numerates lattice sites
$(r = r_1,\,r_2,\,r_3,\,...,r_N)$.
These commutation relations
 imply the independence from
relative arrangement of sites and, therefore,
 essentially restrict the
 domain of values of $q$. In [3] was shown that (1.1)
 requares
 $q=\pm 1$.

 In this paper  we suggest a definition of lattice
  $(l,q)$-deformed boson and Grassmann  fields for
  arbitrary  $q$ and $l$  ($q^2\neq -1$ in the
  last case).
Here $l$ is a deformation parameter for the commutation
relations of "values" of these fields in two arbitrary lattice
sites and $q$ is a deformation parameter for $n$-component
$q$-boson or  $q$-Grassmann variable.

   This paper is organized as follows.
   In section 2 the commutation relations for lattice
  $(l,q)$-deformed fields are
  defined. In section 3 we determine  differential calculus for
  these fields. In section 4 we define the $SL_q(2,C)$-invariant
  extention of Berezin integration for the lattice $(l,q)$-Grassmann
 field.  We show that Gaussian functional integral for the lattice
 $(l,q)$-Grassmann field is expressed in terms of  $(l,q)$-deformed
counterpart of the Pfaffian.

\bigskip

\centerline{\bf 2. Lattice $(l,q)$-deformed fields }

\bigskip

In this section we define  the lattice
 $(l,q)$- deformed boson and Grassmann fields. Let us briefly
 recall the known facts about quantum hyperplanes  $A^{n|0}_q$ and
$A^{0|n}_q$ [12, 13], which for simplicity  we will denote as quantum
vector spaces $\V$ and $\K$ correspondingly.  Define  coordinate
functions for $\V$ as $v^\a$  and correspondingly
$\ka$ for $\K\, (\a=1,\,2,\,...n)$ .
These coordinate functions satisfy the commutation relations which can
write by means of the linear map:
$$
\R\,:\,\, \V\ty\V\longrightarrow \V\ty\V,\,\, \R(v^\a\ty
v^\b)\equiv \R^{\,\a\b}_{\,\g\r}\,v^\g \ty v^\r =q v^\a \ty v^\b
,\tag2.1 $$
 $$
\R\,:\,\, \K\ty\K\longrightarrow \K\ty\K,\,\, \R(\k^\a\ty \k^\b)\equiv
\R^{\,\a\b}_{\,\g\r}\,\kg \ty \kro =-\frac 1q \ka \ty \kb ,\tag2.2
$$
where $\ty$ denotes the tensor product of the quantum spaces,
deformation parameter $q$
in general case is a complex number and the linear map $\R$ is
determined by the symmetric $\R$-matrix of size $n^2 \times n^2$ for
the quantum matrix group $\s$ (q-deformed $SL(n,C)$) [11-13].  The
explicit form of $\R$ is given in [13] , where we see that
$$
\R^{\,\a\b}_{\,\g\r}\,=\R_{\,\a\b}^{\,\g\r},\,\, \R^{\,\a\b}_{\,\g\r}
= \d^\a_\r\,\d^\b_\g (1+(q-1)\d^{\a\b}) + \d^\a_\g\,\d^\b_\r
(q-q^{-1})\Theta^{\a\b}\,,\tag2.3a
$$
$$
\Theta^{\a\b}=\cases
1,&\text {if $\a<\b$}\\ 0,&\text {if $\a\ge \b$}\endcases\,.\tag2.3b
$$
This matrix has eigenvalues $q$ and $-q^{-1}$ and satisfies the
quantum  Yang-Baxter equation [11-13]
$$
\R_{\,12}\,\R_{\,23}\,\R_{\,12}=\R_{\,23}\,\R_{\,12}\,\R_{\,23}\,,\tag2.4a
$$
and the Hecke equation
$$
\R^2-(q-q^{-1})\R-\I=(\I+q\R)(\I-\frac 1q\R)=0, \tag2.4b
$$
where
$
\I=I^{\,\a\b}_{\g\r}=\d^{\a}_{\g}\,\d^{\b}_{\r} $ is unit matrix.
For the $\R$-matrix we can write the following representation
$$
\R=q\,\csq -q^{-1}\caq,
$$
where matrices
$$
\caq  ={\frac {q^2}{1+q^{2}}}(\I-\frac 1q\R), \tag2.5a
 $$
 $$
 \csq  ={\frac 1{1+q^2}}( \I+q\R)  \tag2.5b
$$
are orthonormal projections (relative to the eigenvalues
 $-q^{-1}$ and  $q$ respectively)
 $$
\caq \cdot \csq=\csq \cdot \caq =0,\quad \caq^2=\caq, \quad
\csq^2=\csq.
$$
These projections are the quantum analogs of the
classical $(q=1)$ antisymmetrizer and  symmetrizer for tensors with
two indeces.

Using these matrices,  we can write  commutation
relations (2.1), (2.2) in the form
$$
{{\ca}^{\,\a\b}}_{q\ \g\r}\,v^{\g}\ty v^{\r}=0, \tag2.6
$$
$$
{\cs}^{\,\a\b}_{q\ \g\r}\,\k^{\g}\ty \k^{\r}=0. \tag2.7a
$$
Let us write  last relation  in components  ($q^2\neq -1 $)
$$
\ka\ty\kb=-q\kb\ty\ka\quad \text{for}\,\, \a >\b,\quad
(\ka)^2=0.\tag2.7b
$$
Since at $q=1$ it gives
the commutation relations for the $n$-component classical Grassmann
variable $z^\a$ $$ (z^\a)^2=0,\quad z^\a\,{z^\b}+{z^\b}\,z^\a=0,
 $$
we will consider $\k$ as an $n$-component $q$-Grassmann variable.

Commutation relations (2.6) and (2.7) are
invariant with respect to transformation of $\k$ and $v$ by the quantum
matrix $\A\in\s$ [13]\
$$
\vtt^{\a}=\A^{\a}_{\b}\ty v^{\b},\quad
\kt^{\a}=\A^{\a}_{\b}\ty\k^{\b}, \tag2.8
$$
where matrix elements $\A^{\a}_{\b}$  commute with components
 $\k^{\a}$ and $\va$ and belong to the associative algebra of functions
on the quantum group  $\sq $, which we will denote as $\caa$ or
 $Fun_q(SL\,(n))$ [13]. This algebra is Hopf algebra.  It  means that
in this algebra the following maps are defined:

a) comultiplication $\bigtriangleup$
$$
\caa {\overset\bigtriangleup\to\longrightarrow }\caa\ty\caa\quad:
\qquad \bigtriangleup(A^{\,\a}_{\,\b})=
A^{\,\a}_{\,\g}\ty A^{\,\g}_{\,\b}, $$
where symbol $\ty$ denote the tensor product of quantum space,

b) counit $\e$
$$
\caa {\overset\e\to\longrightarrow } \bold C \quad:\qquad
\e(A^{\,\a}_{\,\b})=\d^{\,\a}_{\,\b},
$$
where $\bold C$ is a complex numbers,

c) antipod $i$
$$
\caa  \overset i\to\longrightarrow \caa \quad:\qquad
i(A^{\,\a}_{\,\b})=
(-q)^{\a-\b}{\tilde A}^{\,\b}_{\,\a},
$$
where ${\tilde A}^{\,\b}_{\,\a} $ is  quantum minor for matrix
element $ A^{\,\a}_{\,\b} $.  For the quantum matrix $\A\in\s$
we have
$$
{\tilde A}^{\a}_{\b}= \sum_{\sigma \in S_{n-1}}
 (-q)^{l(\sigma)}{ A}^{\,1}_{\,\si1 }...{ A}^{\,\a-1}_{\,\si{\a-1}}
{A}^{\,\a+1}_{\,\si{\a+1}}
...\,{ A}^{\,n}_{\,\si{n}}
$$
and
$i(\A)\cdot\A=\A \cdot i(\A)=\hat {\bold 1}$ [13], where "$\cdot$"
denotes the matrix multiplication.
Here the sum extends over the symmetric group $S_{n-1}$, $l(\sigma)$
is the length (the number of inversion) of the substitution
$\sigma= (\si1,\,...,\, \si{\a-1},$\ $
\si{\a+1},\,...,\,\si{n}) =
\sigma (1,\,...,\,\b-1,\b+1,\,...,\,n)$

d) multiplication map $m$ $$
\caa \ty \caa \overset m\to\longrightarrow \caa \quad:\qquad
m(A^{\,\a}_{\,\b} \ty A^{\,\g}_{\,\r})=A^{\,\a}_{\,b}\,
A^{\,\g}_{\,\r}, $$

e) matrix elements  $A^{\,\a}_{\,\b}$ satisfy commutation relations
$$
{\R^{\,\a\b}_{\,\g\r}}\,A^{\,\g}_{\,\mu}\, A^{\,\r}_{\,\nu}=
A^{\,\a}_{\,\g}\, A^{\,\b}_{\,\r}\,{\R^{\,\g\r}_{\,\mu\nu}}.\tag2.9
$$

These maps permit to define the following operations with the quantum
matrix $\A\in\s$:

a)  $\bigtriangleup(\A)=\A\my\A$  defines the rule of multiplication
of quantum matrices and matrix elements of  $\A\my\A$ have form
$A^{\,\a}_{\,\g}\ty A^{\,\g}_{\,\b}$
($\my $ denote the tensor
product of  quantum spases together with  usual matrix
multiplication),

b) $\e(\A)=\hat {\bold 1} $ defines a unit matrix in $\s$,
$\hat {\bold 1}=\d^{\,\a}_{\,\b}$,

c) $i(\A)=\A^{-1}$ defines the inverse matrix,

d)  $m(\A\ty\A)=\A\tt\A$ defines usual tensor multiplication for
quantum matrices and matrix elements of $ \A\tt\A$
have form $A^{\,\a}_{\,\b}\, A^{\,\g}_{\,\r}$.

e)  the quantum determinant of the quantum matrix $\A$\ is determined
by relation
$$\aligned
(A^{\,\a_{1}}_{\,\b_{1}}\,A^{\,\a_{2}}_{\,\b_{2}}\,A^{\,\a_{3}}_{\,\b_{3}}
\cdot\cdot\cdot
A^{\,\a_{n}}_{\,\b_{n}})&\ty(\k^{\b_{1}}\ty\k^{\b_{2}}\ty
\k^{\b_{3}}\ty\cdot\cdot\cdot \ty\k^{\b_{n}})=\\
&det_q \A\,
(\k^{\a_{1}}\ty\k^{\a_{2}}\ty
\k^{\a_{3}}\ty\cdot\cdot\cdot \ty\k^{\a_{n}})\endaligned \tag2.10
 $$
and satisfies property
$$ \bigtriangleup(det_q\A)=det_q(\A\my \A)=det_q(\A)\ty
det_q(\A)\,.
$$
 For $SL_q(n)$
$$
det_q \A = \sum_{\sigma \in S_{n}}
 (-q)^{l(\sigma)}{ A}^{\,1}_{\,\si1 }...\, { A}^{\,n}_{\,\si{n}}=1.
\tag2.11
$$
Using (2.9), it is not hard to show that $q$-determinant commutes with
matrix elements $A^{\,\a}_{\,\b}$ [13].

Now let us consider the definition of the
$(l,q)$-deformed  lattice
fields. At first define  the lattice $n$-component
$(l,q)$-Grassmann field. For this  consider
d-dimensional hypercubic lattice. The lattice  sites are numerated by
radius vector $ r=(p_1,\,p_2,\,.\,.\,.\, ,\,p_d )$, where
$\{p_i\}$ are  integer numbers and  a full number of sites are  $N$ .
For convenience the lattice constant $a$ is fixed to be equal to unit.
Let $\pa{r^{\,\prime}}\in \P_{r^{\,\prime}}$ be the
$n$-component $q$-Grassmann variable $\p^\a$ on the site
$r^{\,\prime}$, which  we will consider as  a "value" of the
$n$-component $(l,q)$-Grassmann field $\pa{r}$ on site $r^{\,\prime}$.
 The definition of the  lattice $(l,q)$-Grassmann
  field requires
 information about the commutation relations of such "values" in
 two arbitrary lattice sites $r_i=(m_1,\,m_2,\,.\,.\,.\, ,\,m_d )$
and $r_j=(k_1,\,k_2,\,.\,.\,.\, ,\,k_d )$.   In
this paper we will suppose the ``quasi''-one-dimensional ordering for the
lattice sites:  $r_j>r_i$, if $k_1-m_1>0$ and the others $k_i-m_i$ are
arbitrary; if $k_1-m_1=0$, then $r_j>r_i$, if $k_2-m_2>0$  and so on,
and  $ r_1\,<\,r_2\,<\,r_3\,<\,...<\,r_N $.

 Let us write these commutation relations in the form
 which
 generalize of (2.7b) in simple way
 $$ \aligned &\pa{r_i} \,\pb{r_i} =
-q \pb{r_i} \, \pa{r_i}\, ,
\quad (\pa{r_i})^2=0,\\ &\pa{r_j} \,\pb{r_j} =
-q \pb{r_j} \, \pa{r_j} \, , \quad (\pa{r_j})^2=0,\\
&\pa{r_j} \, \pa{r_i}
=-l\,\pa{r_i} \, \pa{r_j} \, , \quad\,\, \a
>\b,\,\,r_j>r_i,\endaligned \tag2.12
$$
where ($q^2\neq-1$) and
deformation parameter $l$ in
general case is also a complex number.  At $l=q=1$  we obtain  usual
lattice Grassmann field.

It is not hard to note that
commutation relations (2.12) are covariant with respect to
transformation
$$ {\tilde \psi}^\a_i= A^\a_\b \ty \pb{i} \,, \quad \A
\in \sa\,,
$$ where   $\sa$ is the
quantum matrix group determining the internal properties of the field
$\pa{r}$ and
$$
\tilde \psi^\a_i=L^k_i \ty \pa k\,, \quad \L\in \scl ,
$$
where the space quantum matrix group $\scl$  determines
the commutation properties of the  lattice $q$-Grassmann field  $\pa
r$ in different lattice sites.

Let us denote the quantum vector space which is generated by
$nN$-compo-\-nent $q$-Grassmann vector $\pa r$
as
$\Psi$
($N$ is the  number of sites on lattice).
Commutation relations (2.12) require that
we must consider  quantum vector space $\Psi$ as the Hecke
sum [14,15]  relative to  some matrix $\Q$ which  is defined by
 (2.12)
$$
\Psi=\P_{r_1}{\oplus}_\Q \P_{r_2}{\oplus}_\Q\cdot\cdot\cdot{\oplus}_\Q
\P_{r_{N-1}}{\oplus}_\Q \P_{r_N}\,. \tag 2.13
$$
For definition of   Hecke sum (2.13) and its $\CR$-matrix
 let us consider
 linear map $\Q\,:\,\,\P_{r_i}\ty\P_{r_j}
\longrightarrow \P_{r_j}\ty\P_{r_i}$ \
$$\aligned
&\Q(\pa {r_i}\ty \pb {r_j})\equiv \Q^{\,\a\b}_{\,\g\r}\,\pg {r_i}\ty \pro
{r_j}=-\frac 1l \pa {r_j}\ty \pb {r_i}\,,\\
&\Q^{-1}(\pa {r_j}\ty \pb {r_i})\equiv
{(\Q^{-1})}^{\,\a\b}_{\,\g\r}\,\pg {r_j}\ty \pro {r_i}=-l\, \pa
{r_i}\ty \pb {r_j} \,.\endaligned\tag2.14
$$
where $r_i<r_j\,\,(i<j) $ .
In our case $\Q$ is the matrix of size $n^2\times n^2$
$$
\Q^{\,\a\b}_{\,\g\r} = \d^\a_\g\,\d^\b_\r\, \d^{\a\b} +
q^{-1}\d^\a_\r\,\d^\b_\g \Theta^{\a\b}+ q\d^\a_\r\,\d^\b_\g
\Theta^{\b\a}\,,\tag2.15
$$
and $\Theta^{\a\b}$ is defined in (2.3b). Note that $\Q^{-1}=\Q$.  Let
us show explicit  form of the $Q$-matrix for $n=2$:
 $$
\Q=\pmatrix 1 & 0 &0 &0 \\
\vspace{3truemm} 0& 0 & q & 0 \\
\vspace{3truemm} 0 & q^{-1} & 0 & 0 \\
\vspace{3truemm} 0 & 0 & 0 & 1 \endpmatrix .
$$ By analogy with (2.2), using the matrices $\R$ and $\Q$, we can
define the $\CR$-matrix as linear map $ \CR\,:\,\,
\Psi\ty\Psi\longrightarrow \Psi\ty\Psi$:
$$
\CR(\pa {i}\ty \pb {j})\equiv \CR^{\,\a\b\,k\,l}_{\,\g\r\,i\,j}\,\pg
{k}\ty\pro {l} = - \CI^{\,\a\b\,k\,l}_{\,\g\r\,i\,j}\,\pg {k}\ty\pro
{l}, \tag2.16
$$ where
 here and later on  the values of a latin indeces  goes over all
values of the variable $r$  $(r=r_1,\,r_2,\,r_3,\,...,r_N )$ and
the $\CR$-matrix has the  block structure:
$$ \aligned \CR^{\,\a\b\,k\,l}_{\,\g\r\,i\,j}=&
\R^{\,\a\b}_{\,\g\r}\,\d^k_i\, \d^l_j\,\d^{k\,l}+
(l-l^{-1}) I^{\,\a\b}_{\g\r}\,\d^k_i\,\d^l_j\,
\Theta^{k\,l}+\\ &\Q^{\,\a\b}_{\,\g\r}\,\d^k_j\,\d^l_i\,\Theta^{k\,l}
+\Q^{\,\a\b}_{\,\g\r}\,\d^k_j\,\d^l_i\,\Theta^{l\,k} \,.\endaligned
\tag2.17
$$
$\CI$ is the diagonal matrix
$$
\CI^{\,\a\b\,k\,l}_{\,\g\r\,i\,j}=
\frac 1q \d^\a_\g\,\d^\b_\r\, \d^k_i\, \d^l_j\,\d_{i\,j}-
\frac 1l \d^\a_\g\,\d^\b_\r\, \d^k_i\, \d^l_j\,(\d_{i\,j}-1),
\tag2.18
$$ \
$$
\aligned
&\CI^{\,\a\b\,k\,l}_{\,\g\r\,i\,i}\,\pg {k}\ty\pro {l}=\frac 1q
\,\pa {i}\ty \pb{i},\\
&\CI^{\,\a\b\,k\,l}_{\,\g\r\,i\,j}\,\pg {k}\ty\pro {l}=\frac 1l
\,\pa {i}\ty \pb{j},\quad i\neq j\,.
\endaligned
$$
In (2.16) for  product $\Psi\ty\Psi$ we imply the following
ordering:
$$ \aligned &\Psi\ty\Psi=
\P_{r_1}\ty\P_{r_1}{\oplus}\P_{r_1}\ty\P_{r_2}{\oplus}
\P_{r_2}\ty\P_{r_1}{\oplus}\P_{r_2}\ty\P_{r_2}{\oplus}
\P_{r_1}\ty\P_{r_3}\\ &{\oplus}\P_{r_3}\ty\P_{r_1}{\oplus}
\P_{r_2}\ty\P_{r_3}{\oplus}\P_{r_3}\ty\P_{r_2}{\oplus}
\P_{r_3}\ty\P_{r_3}{\oplus}\cdot\cdot\cdot{\oplus}
\P_{r_N}\ty\P_{r_N}.
 \endaligned\tag2.19
 $$
At such ordering the $\CR$-matrix has  block structure and, for
example, for arbitrary two sites $r_i$ and $r_j\,\,(r_i<r_j)$
the block $\tilde {\Cal R}_{(i\,j)}$ of the $\CR$-matrix acting
on the sum
$\P_{r_i}\ty\P_{r_i}{\oplus}\P_{r_i}\ty\P_{r_j}{\oplus}
\P_{r_j}\ty\P_{r_i}{\oplus}\P_{r_j}\ty\P_{r_j}$
has the following form:
$$
\tilde {\Cal R}_{(i\,j)}=\pmatrix \R & 0 &0
&0 \\
\vspace{3truemm}
0& (l-l^{-1})\I & \Q & 0 \\
\vspace{3truemm}
0 & \Q & 0 & 0 \\
\vspace{3truemm}
0 & 0 & 0 & \R \endpmatrix . \tag2.20
$$

It is not hard to check that the $\CR$-matrix satisfies the quantum
Yang-Baxter equation
$$
\CR_{\,12}\,\CR_{\,23}\,\CR_{\,12}=\CR_{\,23}\,\CR_{\,12}\,\CR_{\,23}
\,,\tag2.21
$$
and the Hecke equation
$$
\CR^2-(\CJ-\CI)\CR-\J=(\J+\CJ\cdot\CR)\cdot(\J-\CI\cdot\CR)=0,
\tag2.22
$$ where
$\J=J^{\,\a\b\,k\,l}_{\,\g\r\,i\,j}=
I^{\,\a\b}_{\g\r}\,I^{\,k\,l}_{i\,j}$ is unit matrix,
 $\CJ$ is the diagonal matrix
$$
\CJ^{\,\a\b\,k\,l}_{\,\g\r\,i\,j}=
q\, \d^\a_\g\,\d^\b_\r\, \d^k_i\, \d^l_j\,\d_{i\,j}-
l\, \d^\a_\g\,\d^\b_\r\, \d^k_i\, \d^l_j\,(\d_{i\,j}-1),
\tag2.23
$$\
$$
\aligned
&\CJ^{\,\a\b\,k\,l}_{\,\g\r\,i\,i}\,\pg {k}\ty\pro {l}=q
\,\pa {i}\ty \pb{i},\\
&\CJ^{\,\a\b\,k\,l}_{\,\g\r\,i\,j}\,\pg {k}\ty\pro {l}=l
\,\pa {i}\ty \pb{j},\quad i\neq j\,.
\endaligned
$$
and $\J=\CI\cdot\CJ$.
By analogy with (2.5a) and (2.5b), from (2.22) we can define
  the quantum $(l,q)$-antisymmetrizer
$$
\baq  =(\J-\CI\cdot \CR) \tag2.24
 $$
 and  the quantum $(l,q)$-symmetrizer
 $$
 \bsq  =(\J+\CJ\cdot\CR),  \tag2.25
$$
which are orthogonal projections
 $$\aligned
&\baq \cdot \bsq=\bsq \cdot \baq =0,\\
&(\baq)^2=\CI^2\cdot(\CI+\CJ)\cdot\baq, \\
&(\bsq)^2=\CJ^2\cdot(\CI+\CJ)\cdot\bsq.\endaligned
\tag 2.26
$$
The quantum
symmetrizer $\bsq$ permits to define the commutation relations for the
  lattice $n$-component $(l,q)$-Grassmann field $\pa r$ $(q^2\neq -1)$
which include  (2.12)
$$ \bsq\cdot (\psi\ty\psi)=0,\,\,\text{or}\,\,
\pa {i}\ty\pb {j}+
(\CJ\cdot\CR)^{\,\a\b\,k\,l}_{\,\g\r\,i\,j}\,\pg {k}\ty\pro {l}=0\,,
\tag2.27
$$
or in components
$$
\aligned &\pa{r_i} \ty\pb{r_i} = -q\, \pb{r_i} \ty \pa{r_i}\, , \quad
(\pa{r_i})^2=0,\\
&\pa{r_j} \ty\pa{r_i} = -l\, \pa{r_i} \ty \pa{r_j}\, ,\\
&\pa{r_j} \ty\pb{r_i} = -l\,q \, \pb{r_i} \ty \pa{r_j}\, ,\\
&\pb{r_j} \ty\pa{r_i} = -\frac lq  \, \pa{r_i} \ty \pb{r_j}\, ,     \\
\endaligned     \tag2.28
$$
where $r_i<r_j$ and $\a>\b$.

  The quantum antisymmetrizer $\baq$ determines the commutation
relations for the lattice $n$-component $(l,q)$-boson field
$\vf=\{\,\fa r\,\}$
$$
\baq\cdot (\vf\ty\vf)=0\,\,\text{or}\,\,
\fa {i}\ty\fb {j}
-(\CI\cdot \CR)^{\,\a\b\,k\,l}_{\,\g\r\,i\,j}\,\fg {k}\ty\fro
{l}=0\,,\tag 2.29
$$
or in components
$$ \aligned
&\fa{r_i} \ty\fb{r_i} = \frac 1q  \fb{r_i} \ty \fa{r_i}\, , \\
&\fa{r_j} \ty\fa{r_i} = \frac 1l  \fa{r_i} \ty \fa{r_j}\, , \\
&\fa{r_j} \ty\fb{r_i} = \frac q{l}  \fb{r_i} \ty \fa{r_j}\, , \\
&\fb{r_j} \ty\fa{r_i} = \frac 1{l\,q} \,\fa{r_i} \ty \fb{r_j}\, .
\endaligned     \tag2.30
$$
Here again
 $r_i<r_j$ and $\a>\b$.
In this case  $\Q$-matrix (2.15)  is acting by the
following way:
$$\aligned
&\Q(\fa {r_i}\ty \fb {r_j})\equiv \Q^{\,\a\b}_{\,\g\r}\,\fg {r_i}\ty \fro
{r_j}=l\, \fa {r_j}\ty \fb {r_i}\,,\\
&\Q^{-1}(\fa {r_j}\ty \fb {r_i})\equiv
{(\Q^{-1})}^{\,\a\b}_{\,\g\r}\,\fg {r_j}\ty \fro {r_i}=\frac 1l  \,
\fa {r_i}\ty \fb {r_j} \,.\endaligned
$$
Note that $\CR$-matrix (2.24) in which  we   substitute  the
transpose matrix $\Q^t$ insead of $\Q$-matrix (2.15) also satisfies
the quantum Yang-Baxter equation (2.21) and the Hecke equation (2.22).
In this case commutation relations (2.28) for
the lattice $n$-component $(l,q)$-Grassmann field  and (2.30) for
the lattice $n$-component $(l,q)$-boson field  take on the following form
$$
\aligned &\pa{r_i} \ty\pb{r_i} = -q\, \pb{r_i} \ty \pa{r_i}\, , \\
&\pa{r_j} \ty\pa{r_i} = -l\, \pa{r_i} \ty \pa{r_j}\, ,\\
&\pa{r_j} \ty\pb{r_i} = -\frac lq \, \pb{r_i} \ty \pa{r_j}\, ,\\
&\pb{r_j} \ty\pa{r_i} = - l\,q  \, \pa{r_i} \ty \pb{r_j}\, ,     \\
\endaligned
$$
and
$$ \aligned
&\fa{r_i} \ty\fb{r_i} = \frac 1q  \fb{r_i} \ty \fa{r_i}\, , \\
&\fa{r_j} \ty\fa{r_i} = \frac 1l  \fa{r_i} \ty \fa{r_j}\, , \\
&\fa{r_j} \ty\fb{r_i} = \frac 1{l\,q}  \fb{r_i} \ty \fa{r_j}\, , \\
&\fb{r_j} \ty\fa{r_i} = \frac q{l} \,\fa{r_i} \ty \fb{r_j}\, ,
\endaligned
$$
where $r_i<r_j$ and $\a>\b$.

Now let us consider a general linear transformation of
 the lattice field $\pa r$   by the
 matrix $\F$
$$
\Psi  \overset\d\to\longrightarrow
\F\my\Psi:\, \tp = \d(\p)=\F\my\p, \,\,  \tp^\a_{i} =\d(\pa
{r})=F^{\a\,k}_{\b\, i}\ty\pb{k},
\tag2.31
$$
and suppose that the  matrix elements
$F^{\a\,k}_{\b\, i}$ belong to the associative algebra of functions
$\CF$ for which are defined the following operations:

a) comultiplication $\bigtriangleup$
$$
\CF {\overset\bigtriangleup\to\longrightarrow }\CF\ty\CF\quad:
\qquad \bigtriangleup(F^{\a\,k}_{\b\, i})=
F^{\a\,l}_{\r\, i}\ty F^{\r\,k}_{\b\, l},
$$

b) counit $\e$
$$
\CF {\overset\e\to\longrightarrow } \bold C \quad:\qquad
\e(F^{\r\,k}_{\b\, i})=\d^{\,\a}_{\,\b}\,\d^{\,k}_{\,i},
$$
where $\bold C$ is a complex numbers,

c) antipod $i$
$$
\CF  \overset i\to\longrightarrow \CF \quad:\qquad
i(F^{\a\,k}_{\b\, i})=
(-q)^{(\a-\b)+(k-i)}{\tilde F}^{\b\,i}_{\a\,k},
$$
where ${\tilde F}^{\b\,i}_{\a\,k} $ is the quantum minor for matrix
element $ F^{\a\,k}_{\b\, i} $,

d) multiplication map $m$
$$
\CF \ty \CF \overset m\to\longrightarrow \CF \quad:\qquad
m(F^{\a\,k}_{\b\, i} \ty F^{\r\,l}_{\g\, j})=F^{\a\,k}_{\b\, i}\,
F^{\r\,l}_{\g\, j}, $$

e)  matrix elements  $F^{\a\,k}_{\b\, i}$ satisfy
the commutation relations
$$
\CR^{\,\a\b\,k\,l}_{\,\g\r\,i\,j}\,F^{\g\,m}_{\mu\, k}\,
F^{\r\,n}_{\nu\, l}= F^{\a\,k}_{\g\, i}\,F^{\b\,l}_{\r\, j}\,
\CR^{\,\g\r\,m\,n}_{\,\mu\nu\,k\,l}\,.\tag2.32
$$
By analogy with (2.10) one defines $(l,q)$-determinant for matrix
$\F$:
$$\aligned
&(F^{\a_{1}\,k_{1}}_{\b_{1}\,
i_{1}}\cdot\cdot\cdot F^{\a_{n}\,k_{n}}_{\b_{n}\, i_{1}} \,
F^{\,\n_{1}\,l_{1}}_{\g_{1}\,
i_{2}}\cdot\cdot\cdot F^{\,\n_{n}\,l_{n}}_{\g_{n}\, i_{2}}
\cdot\cdot\cdot F^{\,\d_{1}\,m_{1}}_{\m_{1}\, i_{N}}\cdot\cdot\cdot
F^{\,\d_{n}\,m_{N}}_{\m_{n}\, i_{N}})\ty\\
&(\p^{\,\b_{1}}_{k_{1}}\ty\cdot\cdot\cdot\ty \p^{\,\b_{n}}_{k_{n}}\ty
\p^{\,\g_{1}}_{l_{1}}\ty\cdot\cdot\cdot\ty
\p^{\,\g_{n}}_{l_{n}}\ty\cdot\cdot\cdot\ty
\p^{\,\m_{1}}_{m_{1}}\ty\cdot\cdot\cdot\ty \p^{\,\m_{n}}_{m_{n}})=\\
&det_{(l,q)}\F\,
(\p^{\a_{1}}_{i_{1}}\ty\cdot\cdot\cdot\ty \p^{\a_{n}}_{i_{1}}\ty
\p^{\,\n_{1}}_{i_{2}}\ty\cdot\cdot\cdot\ty
\p^{\,\n_{n}}_{i_{2}}\ty\cdot\cdot\cdot\ty
\p^{\d_{1}}_{i_{N}}\ty\cdot\cdot\cdot\ty \p^{\d_{n}}_{i_{N}})\,.
\endaligned\tag2.33
$$
Then  it is not hard show that $\CF$ is  algebra of functions
$Fun_{(l,q)} (GL(n)\times GL(N))$ and
 commutation relations (2.28) and (2.29) are covariant with
respect to transformation (2.31). This algebra contains
the algebras \hfill\  \linebreak $Fun_{q} (GL(n))$
and $Fun_{l} (GL(N))$ as subalgebras.

For futher we restrict ourselves to such transformations (2.31)
which do not depend from the coordinates of lattice sites
$$
F^{\a\,k}_{\g\, i}=A^{\,\a}_{\,\g}\,\d^{\,k}_{\,i}\,.\tag2.34
$$
Let us suppose that $det_q\A=1$. Then from commutation
relations (2.32) and the expression (2.17) for $\CR$ it follows
that $A^{\,\a}_{\,\g}\in Fun_q(SL(n))$.
 It means that $\caa\in\CF$.
Really, taking into account (2.34), from (2.33) and (2.10) it is not
hard to show that $$ det_{(l,q)}\F=(det_q\A)^N=1.\tag2.35 $$

\bigskip

\centerline{\bf 3. Differential calculus for}
\centerline{\bf  the lattice $(l,q)$-deformed fields}

\bigskip

The various approaches to the noncommutative differential geometry
on the quantum spaces have been considered in  [16-20] (for a review
see [20]).  In this section, following the Wess and Zumino approach
[16,17], we consider the differential calculus for
lattice $(l,q)$-deformed fields defined in previous section.
At first
 let us determine differential calculus for the lattice
$(l,q)$-Grassmann field $\p$ satisfying commutation
relations (2.27).

Define a exterior differential $d$ on $\Psi$ satisfying the
usual properties such as
$$
d^2=0 \tag3.1
$$
and the Leibniz rule
$$
d(fg)=(df)g+(-1)^{|f|}f(dg)\,, \tag3.2
$$
where  $|f|$ denotes the parity of $f$  (for example, $|\p|=-1$) and
$f$ and  $g$ are a functions of the $(l,q)$-Grassmann field $\p$ .
Under these functions we understand their formal expansion in power
series of  $\p$, which have the finite number of members because of
the nilpotency condition (2.13)  for  $\pa{r}$.
 We will consider $ \{\psi^{\a}_r\}$  as
the basic elements for $ \P$.  In general case for
derivation of the commutation relations between functions,
  differentials and derivatives it is necessary to consider their
 action on arbitrary function $f(\p)$.
However for $\P$, in consequence of the finite number of
 members of the power series  and the form of commutation
 relations (2.28) we can restrict ourselves to study of their action
  on basic element $ \{\psi^{\a}_r\}$.

  The action  of a
exterior differential on this field is
$$
d(\pa{r}\,\pb{r\prime})=(d\pa{r})\,\pb{r\prime}-\pa{r}\,(d\pb{r\prime})
\,. \tag3.3
$$
Define the derivative operator
$$
\ia k =
\,\dfrac{\partial\ }{\partial\pa{k}}\,,\quad
 \ia{l}\,\pb{k} =\d^{\b}_{\a}\,\d^{l}_{k},\quad
 d=d\pa{k}\,\ia{k}\,, \tag3.4
$$
 Write the commutation relations between the lattice field
 $\p$  and its differential $d\p$ in the following general form
 $$
 \p\ty d\p = \C\cdot (d\p\ty\p)=(\J+\B)\cdot
 (d\p\ty\p)\,, \tag3.5
 $$
 where $\C$ and $\B$ are  some numerical
matrices, which we will find from the consistency condition of (3.5)
with (2.27).

Acting by  exterior differential $d$ on (2.27)
$$
 \bsq\cdot (\p\ty d\p-d\p\ty \p) =0\,,
 $$
 and substituting in this equation  (3.5) we obtain the following
 consistency condition
 $$
 \bsq\cdot\B=0\,. \tag3.6
 $$
 Using the orthogonality property of
  $\baq$ and $\bsq$
  (2.26) we suggest the following anzats for solution this condition
 $$ \B=a_1\,\baq.
 $$
 Requiring that  commutation relation (3.5) at $\a\neq\b,\,i=k$
  coinsides with
 $$
\pa{r}\,d\pb{r}  -{\frac 1 q}\,\mr d\pg{r}\,\pro{r}=0\, , $$
 we can determinate the coefficient $a_1$:
 $a_1=-1$.
 Then
 $$ \aligned
  \C=
  {\CI}\cdot\CR
  \endaligned\tag3.7
  $$
and we obtain
 $$
\pa{i}\,d\pb{j}  -(\CI\cdot\CR)^{\a\b\,k\,l}_{\g\r\,i\, j}
 d\pg{k}\,\pro{l}=0\,.\tag3.8
 $$

Note that the derivatives defined by  relation (3.4)
do not satisfy the Leibniz rule.
Indeed using (3.3) and (3.8), we obtain

$$
\split
&d(\pa{i}\,\pb{j})=
(d\pa{i})\,\pb{j}-\pa{i}\,(d\pb{j})=\\
&d\pg{k}(\ig{k}\,\pa{i})\pb{j}-\pa{i}\,d\pg{k}(\ig{k}\,\pb{j})=\\
&d\pg{k}\lbrack (\ig{k}\,\pa{i})\pb{j}-(\CI\cdot\CR)^{\,\a\d\,k\,
m}_{\,\g\r\,i\, l}\,\pro{m}(\id{l}\,\pb{j})\rbrack =\\
&d\pg{k}\,\ig{k}(\pa{i}\,\pb{j}).
\endsplit
$$

Hence it follows
$$
\ig{k}(\pa{i}\,\pb{j})=(\ig{k}\,\pa{i})\pb{j} -
(\CI\cdot\CR)^{\,\a\d\,k\, m}_{\,\g\r\,i\, l}\,
\pro{m}(\id{l}\,\pb{j}) .\tag3.9
$$

If consider $\ia{i}$ and $\pb{j}$ as operators, acting on some
 function $f(\p)$, then from  (3.9) we obtain commutation relations
 $$
\ig{k}\,\pa{i}+
(\CI\cdot\CR)^{\,\a\b\,k\, l}_{\,\g\r\,i\, j}\,
\pro{l}\,\ib{j}=
\d^{\a}_{\g}\,\d^{k}_{i}.\tag3.10
$$
Let us assume that the  matrix $\D$ determines the commutation
relations between derivatives $\ig{k}$ and differentials $d\pa{j}$
 $$
 \ig{k}\,d\pa{i}-
{\D}^{\,\a\b\,k\, i}_{\,\g\r\,j\, m}\,
d\pro{i}\,\ib{m}=  0.
 \tag3.11
 $$
To find the matrix
$\D$  let us consider (3.11) as the operator identity and  act by
them on $\pd{m}$.  Using (3.8) and (3.10), we obtain
$$
\split
& \ig{k}\,d\pa{j}\,\pd{m} - {\D}^{\,\a\b\,k\, i}_{\,\g\r\,j\, l}\,
d \pro{i}(\ib{l}\pd{m})=\\
&{[(\CI\cdot\CR)^{-1}]}^{\,\a\d\,i\, l}_{\,\mu\nu\,j\, m}\,
(\ig{k}\,\pmm{i})d\pn{l} -{\D}^{\,\a\b\,k\, i}_{\,\g\r\,j\, l}\,
d \pro{i}(\
ib{l}\pd{m})=\\
&{[(\CI\cdot\CR)^{-1}]}^{\,\a\d\,i\, l}_{\,\mu\nu\,j\, m}\,
(\d^\mu_\g\,\d_{i}^k)d\pn{l} -
{\D}^{\,\a\b\,k\, i}_{\,\g\r\,j\,l}\,d \pro{i}(\d_\b^\d\,\d_{m}^l)=\\
&[{[(\CI\cdot\CR)^{-1}]}^{\,\a\d\,k\, i}_{\,\g\r\,j\, l}\,
 - {\D}^{\,\a\b\,k\, i}_{\,\g\r\,j\,l}]\,
d \pro{i}=0.
\endsplit
$$
It implies  $\D=(\CI\cdot\CR)^{-1}$ and we get
$$
d\pg{i}\,\ia{k} -(\CI\cdot \CR)^{\,\a\d\,k\, l}_{\,\g\r\,i\,j}\,
\id{j}(d\pro{l})= 0.  \tag3.12
$$
In order to determine the commutation relations between the
derivatives we use the following anzats
$$
\CG^{\a\b\,k\,l}_{\,\g\r\,i\,j}\,\ia{i}\,\ib{j}=0,\tag3.13
$$
where $\CG$ is some numerical matrix.  Multilying this equation from
right by $\pd{m}$ and commuting $\pd{m}$ through to the left and
requiring that terms linear in the derivatives are canceled, we find
the consistency condition
$$
(\J- \CI\cdot\CR)\cdot\CG=0 \quad {\text or}\quad\baq\cdot
\CG=0.\tag3.14
$$
Using   orthogonality property (2.26)
and again requiring that the commutation relation (3.13) at $\a\neq\b,\,i=k$
  has to  coincide with
 $$
\ia{i}\,\ib{i}  + q\,\R^{\,\g\r}_{\,\a\b}\, \ig{i}\,\iro{i}=0\, ,
 $$
 we  find:
 $$
\CG =\bsq.
 $$

Hence we obtain  commutation relations:
 $$
\aligned
&\ia{r_i} \,\ib{r_i} = - q\,  \ib{r_i} \, \ia{r_i}\, , \\
&\ia{r_j} \,\ia{r_i} = - l\,  \ia{r_i} \,\ia{r_j} \, \\
&\ia{r_j}
\,\ib{r_i}  =-{l\,q}\, \ib{r_i} \,\ia{r_j}  \, \\
&\ib{r_j}
\,\ia{r_j}  =-\frac lq \,\ia{r_i} \,\ib{r_j} \, , \endaligned
\tag3.15
$$
where again $\a>\b$ and $r_j>r_i$ or in covariant form
$$
\ia{k}\,\ib{l}  + (\CJ\cdot\CR)^{\,\g\r\,k\,l}_{\,\a\b\,i\,j}\,
\ig{i}\,\iro{j} =0\,.\tag3.16
$$

Now it is not hard to find  commutation relations between
differentials. For this aim,
acting by
operator $d$ on (3.8)
and using the  relation
$$
d^2(\pa{r}\,\pb{r\prime}) =
d(d\pa{r}\,\pb{r\prime}-\pa{r}\,d\pb{r\prime})=
d(d\pa{r}\,\pb{r\prime})-d\pa{r}\,d\pb{r\prime}=0,
$$
we obtain
 $$
d\pa{i}\,d\pb{j}  -(\CI\cdot\CR)^{\a\b\,k\, l}_{\g\r\,i\, j}
 d\pg{k}\,d\pro{l}=0\,.\tag3.17
 $$
By analogy with the  determination of the
commutation relations
(3.8), (3.10), (3.12), (3.16), (3.17)
for the lattice $n$-component $(l,q)$-Grassmann field
we
can find analogous commutation relations for the lattice $n$-component
$(l,q)$-boson field which satisfies commutation relations (2.29).
 Here we give the final result:
$$
\align
&\fa{i}\,d\fb{
j}  -(\CJ\cdot\CR)^{\a\b\,k\,l}_{\g\r\,i\, j}
 d\fg{k}\,\fro{l}=0\,,\tag3.18\\
&\ig{k}\,\fa{i}-
(\CJ\cdot\CR)^{\,\a\b\,k\, l}_{\,\g\r\,i\, j}\,
\fro{l}\,\ib{j}=
\d^{\a}_{\g}\,\d^{k}_{i}\,,\tag3.19\\
&d\fg{i}\,\ia{k} - (\CJ\cdot\CR)^{\,\a\b\,k\,
l}_{\,\g\r\,i\,j}\, \ib{j}(d\fro{l})= 0\,,  \tag3.20\\
&\ia{k}\,\ib{l}  - (\CI\cdot\CR)^{\,\g\r\,k\,l}_{\,\a\b\,i\,j}\,
\ig{i}\,\iro{j}=0 \,,\tag3.21\\
&d\fa{i}\,d\fb{j}  + (\CJ\cdot\CR)^{\a\b\,k\, l}_{\g\r\,i\, j}
 d\fg{k}\,d\fro{l}=0\,,\tag3.22
\endalign
$$
where we use the notation
$$
\ia k =
\,\dfrac{\partial\ }{\partial\fa{k}}\,,\quad
 \ia{l}\,\fb{k} =\d^{\b}_{\a}\,\d^{l}_{k},\quad
 d=d\fa{k}\,\ia{k}\,. \tag3.23
 $$

\bigskip

\centerline{\bf 4. The elements of integral calculus for}
\centerline{\bf  the lattice $(l,q)$-Grassmann fields}

\bigskip

In  this section we determine the generalization of the Berezin
integration  for the lattice $(l,q)$-Grassmann field defined in
previous sections.

At first let us  define
the $\s$-invariant generalization of the Berezin integration   on
the quantum vector space $\K$.
For this aim we will use  the results  of the paper [3] where
 representation of the partition function of the two-dimensional
Ising model in form of the $q$-Grassmann functional integral ($q=-1$)
was found.
After calculation of this integral for the Ising model with nearest
neighbour interaction   the known  solution of the model was obtained
  in [3].

  Let us define  the volume form  which is invariant
with   respect to  transformation from $\s$.
 For that it is necessary to
introduce  the operation of the exterior multiplication " $\wedge $ ".
  Make this, using the $q$-symmetrizer $\csq $ (2.5b)
  $$ \cs
^{\,\,\a\b}_{\,q\ \g\r}\,d\kg\w d\kr\equiv 0. \tag 4.1
$$
 This definition
permits to determine the commutation relations between the
differentials $d\ka$ relative to the exterior multiplication
$$
d\ka\w d\kb=-q\,\mr d\kg\w d\kr\,.\tag4.2
$$
Using this relation, it is easy to show that the volume form
$d\k^n\w d\k^{(n-1)}\w ...d\k^2\w d\k^1\,$ is invariant with respect
to the transformation $\A\in \s $\ $$ \aligned
&d\kt^n\w d\kt^{(n-1)}\w ...d\kt^2\w d\kt^1=\\ &det_q\A \,
d\k^n\w d\k^{(n-1)}\w ...d\k^2\w d\k^1 = d\k^n\w d\k^{(n-1)}\w ...d\k^2\w
d\k^1.  \endaligned \tag4.3
$$
Taking into account this invariance and requiring that  at $q=1$ we
must obtain  the rules of Berezin inteqration [21] for usual Grassmann
variable, one defines
$$
\int d\k^n\w d\k^{(n-1)}\w ...d\k^2\w d\k^1=0. \tag4.5
$$

Note that the differentials in our  volume form (4.3) are different
from Berezin differentials $(q=1)$ since the latter is transformed by
the reciprocal matrix. We are grateful to D.V.~Volkov for this remark.

In consequence of  commutation relations (2.2) for $\ka$,
it is easy to show, that product
$\k^n\, \k^{(n-1)}\, ...\k^2\, \k^1 $ is also invarint with
respect to transformations from $\s$$$
\kt^1\, \kt^2\, ...\kt^{(n-1)}\, \kt^n =
\k^1\, \k^2\, ...\k^{(n-1)}\, \k^n.
 $$
Using this invariance and (4.3), and again assumig the right limit at
$q=1$, one defines
$$
\int\,d\k^n\w d\k^{(n-1)}\w ...d\k^2\w d\k^1\,
\k^1\, \k^2\, ...\k^{(n-1)}\k^n=1.\tag4.6
$$
Analogous reasoning leads to the definition
$$
\int\,d\k^n\w d\k^{(n-1)}\w ...d\k^2\w d\k^1\,
\k^1\, \k^2\, ...\k^{(n-1)}=0.\tag4.7
$$
The definitions (4.5)--(4.7) are able to write in the following
covariant form
$$
\int\,d\k^n\w d\k^{(n-1)}\w ...d\k^2\w d\k^1\, \k^{\a_1}\, \k^{\a_2}\,
...\k^{\a_{(n-1)}}\,\k^{\a_n}= \e^{\,\a_1\a_2 ...\a_{(n-1)}\a_n}
\,,\tag4.8 $$
Here tensor $\e^{\,\a_1\a_2 ...\a_{(n-1)}\a_n}$ is defined by the
following rules:  $$
\e^{\,1\,2\,3\,...{(n-1)}\,n}=1,\tag4.9 $$
and the remaining components are defined by the coefficients appearing
in the l.h.s of the relation
$$
\k^{\a_1}\, \k^{\a_2}\, ...\k^{\a_{(n-1)}}\,\k^{\a_n}= \e^{\,\a_1\a_2
...\a_{(n-1)}\a_n}\,\k^1\, \k^2\, ...\k^{(n-1)}\k^n\tag4.10
$$
after reordering of its the l.h.s. to
the r.h.s. by means of the commutation relations (2.2) for every
specific set of the values of indices.

Note that in [23,24], using differential and integral calculi
on the quantum plane which are
invariant with respect to quantum
inhomogeneous Euclidean group
 $E(2)_q$, the holomorphic representation
of $q$-deformed path integral for the quantum mechanical evolution
operator kernel of q-oscillator was constructed .

For the calculation of the Gaussian integral over the $n$-component
 $q$-Grassmann variable  we will use the following definition
of  $q$-deformed Pfaffian (which is a generalization
of  definition of usual
Pfaffian [25]
for case of $q$-deformed
commutation relations).

   Consider the quadratic form
$$
w=\,
\sum_{\a<\b}\,a^{\a\b}\,\k^\a\,\k^\b=\,
{\frac 12}\sum_{\a,\b}\,{\hat a}^{\a\b}\,\k^\a\,\k^\b\,, \tag4.11
 $$
where  matrix elements ${\hat a}^{\a\b}$
commute between themselves and $q$-antisymmetric matrix $\hat a$ has
form
$$
\hat a =\bmatrix
0 & {\hat a}^{12}&{\hat a}^{13}&...&{\hat a}^{1n}\\
\vspace{3truemm}
-q^{-1}{\hat a}^{12}&0&{\hat a}^{23}&...&{\hat a}^{2n}\\
\vspace{3truemm}
-q^{-1}{\hat a}^{13}&-q^{-1}{\hat a}^{23}&0&...&{\hat a}^{3n}\\
\vspace{3truemm}
...&...&...&...&...\\
\vspace{3truemm}
-q^{-1}{\hat a}^{1n}&-q^{-1}{\hat a}^{2n}&-q^{-1}{\hat
a}^{3n}&...&0 \endbmatrix\, . \tag4.11a
$$
Here the matrix elements arranged down the diagonal are determined
by means of relations (2.2).
Define $ q$-Pfaffian of matrix $\hat a$ through
$$
\frac 1{({\frac n2})!}\, w^{\frac n2}= Pf_q (\hat a)\,\k^1\, \k^2\,
...\k^{(n-1)}\k^n\, .  \tag4.12 $$
Hence,  taking into account
commutation relations (2.2), we get,  for example, at $n=4$ $$
Pf_q(\hat a)={\frac 12}(1+q^4){\hat a}^{12}\,{\hat a}^{34}- {\frac
1{2}}q(1+q^2){\hat a}^{13}\,{\hat a}^{24}+q^2 {\hat a}^{14}\, {\hat a}^{23}\,.
$$
Using (4.12), it is not hard to show   that
$$
\int\,d\k^n\w d\k^{(n-1)}\w ...d\k^2\w d\k^1\,
\exp \biggl\{{\frac
12} \sum_{\a,\,\b}\,{\hat a}^{\a\b}\,\k^\a\,\k^\b \biggr \}
= Pf_q(\hat a).\tag4.13
$$
Now let us define the integral calculus for
for the lattice $n$-component
$(l,q)$-Grassmann
 field $\p=\{\p^{\a}_r\}$ which satisfy commutation
relations  (2.28). Recall that the general linear
transformation for this field is transformation by the quantum matrix
 $\F\in Fun_{(l,q)}( GL(N) \times GL(n))$.  Following to (2.34),
 we restrict ourselves to transformations $\F$ which do not depend
 from coordinates of lattice sites and belong to
  $Fun_q(SL(n))$ that it means  $det_q(\F)=1$ (2.35).

  By
analogy with (4.1), (4.2), let us define the commutation relations
 between the differentials $d\pa{r}$ relative to the exterior
 multiplication
 $$ (\bsq) ^{\,\,\a\b\,k\,l}_{\,q\
\g\r\,i\,j}\,d\pg{k}\w d\pro{l}=
   \CJ^{\,\a\b\,k\,l}_{\,\g\r\,i\,j}\,\pa {i}\w\pb {j}
+(\CJ\cdot\CR)^{\a\b\,k\, l}_{\g\r\,i\, j} d\pg{k}\w d\pro{l}\equiv 0.
\tag 4.14 $$ or in components
$$ \aligned &d\pa{r_i} \w d\pb{r_i} = - q\,  d\pb{r_i} \w d\pa{r_i}\,
, \\ &d\pa{r_j} \w d\pa{r_i} = - l\,  d\pa{r_i} \w d\pa{r_j} \, ,\\
&d\pa{r_j} \w d\pb{r_i}  =- {l\,q}\, d\pb{r_i} \w d\pa{r_j}  \,,\\
&d\pa{r_i} \w d\pb{r_j}  =- \frac ql d\pb{r_j} \w d\pa{r_i} \, ,
\endaligned \tag4.15
$$
where again $\a>\b$ and $r_j>r_i$.

The exterior multiplication permits to
define the $SL_q(n,C)$-invariant volume form $\cD\p$ with respect to the
transformation by the quantum matrix $\F$ (2.34):
$$ \aligned &\cD\p=
d\p^n_{r_N}\w d\p^{(n-1)}_{r_N}\w ...d\p^2_{r_N}\w d\p^1_{r_N}\w ...\\
&d\p^n_{r_2}\w d\p^{(n-1)}_{r_2}\w ...d\p^2_{r_2}\w d\p^1_{r_2}
d\p^n_{r_1}\w d\p^{(n-1)}_{r_1}\w ...d\p^2_{r_1}\w
d\p^1_{r_1}\,,\endaligned \tag4.16
$$
Then we can define the $SL_q(n,C)$-invariant generalization of the
Berezin integration for the  lattice $(l,q)$-deformed field $\p^\a_r$
by the following relation $$\aligned &\int\cD\p\, \p^{\a_1}_{i_1}\,
\p^{\a_2}_{i_1}\, ...\,\p^{\a_{n-1}}_{i_1}\, \p^{\a_n}_{i_1} ...
\,\p^{\d_1}_{i_N}\, \p^{\d_2}_{i_N}\,
...\,\p^{\d_{n-1}}_{i_N}\,\p^{\d_n}_{i_N}= \\
&\e^{\,\a_1\a_2 ...\a_{n-1}\a_n ...\,\d_1\d_2 ...\,\d_{n-1}
\d_n }_{\,i_1\,i_2 ...\,i_{N-1}\,i_N} \,, \endaligned \tag4.17
$$
where tensor $\hat \e$ is defined by  rules:
$$
\e^{\,1\,2 ...\,{(n-1)}\,n ...\,1\,2 ...\,{(n-1)}
\,n }_{\,r_1\,r_2 ...\,r_{N-1}\,r_N}=1\,,
$$
and the remaining components are defined by the coefficients appearing in
the l.h.s. of the relation

$$  \aligned
&\p^{\a_1}_{i_1}\, \p^{\a_2}_{i_1}\, ...\p^{\a_{n-1}}_{i_1}\,
\p^{\a_n}_{i_1} ...
\,\p^{\d_1}_{i_N}\, \p^{\d_2}_{i_N}\,
...\,\p^{\d_{n-1}}_{i_N}\,\p^{\d_n}_{i_N}= \\
&\e^{\,\a_1\a_2 ...\a_{n-1}\a_n ...\d_1\d_2 ...\d_{n-1}
\d_n }_{\,i_1\,i_2 ...\,i_{N-1}\,i_N} \,
\p^{\,1}_{r_1}\, \p^{\,2}_{r_1}\, ...\p^{\,{n-1}}_{r_1}\,
\p^{\,n}_{r_1} ...\,
\p^{\,1}_{r_N}\, \p^{\,2}_{r_N}\,
...\,\p^{\,{n-1}}_{r_N}\,\p^{\,n}_{r_N}\endaligned
$$
after reordering of its the l.h.s. to
the r.h.s. by means of the commutation relations (2.28) for every
specific set of the values of indeces.

For the calculation of the Gaussian integral over the lattice
 $(l,q)$-Grassmann field  we will use the following definition
of  $(l,q)$-Pfaffian
which is similar  to (4.11).

   Consider the quadratic form
$$
w=\,
\sum_{\a<\b\,i<k}\,b^{\a\b}_{i\,k}\,\p^\a_i\,\p^\b_k=\,
{\frac 12}\sum_{\a,\b,i,k}\,{\hat
 b}^{\a\b}_{i\,k}\,\p^\a_i\,\p^\b_k\,, \tag4.18
 $$
 where
  matrix elements ${\hat b}^{\a\b}_{i\,k}$ commute between
themselves and $(l,q)$-antisymmetric matrix $\hat b$ is determined
by means of the commutation relations (2.29), (2.30) and has the form
which is similar to (4.11a).  Define a $(l,q)$-Pfaffian of matrix $\hat
b$ through $$ \frac 1{({\frac {nN}2})!}\, w^{\frac{nN}2}= Pf_{(l,q)}
(\hat b)\, \p^{\,1}_{r_1}\, \p^{\,2}_{r_1}\, ...\p^{\,{n-1}}_{r_1}\,
\p^{\,n}_{r_1} ...\, \p^{\,1}_{r_N}\, \p^{\,2}_{r_N}\,
...\,\p^{\,{n-1}}_{r_N}\,\p^{\,n}_{r_N} \tag4.19
$$
Using (4.17) and
(4.19), it is not hard to show   that $$ \int\,\cD\p \exp
\biggl\{{\frac 12} \sum_{\a,\b,i,k}\,{\hat
b}^{\a\b}_{i\,k}\,\p^\a_i\,\p^\b_k \biggr \} = Pf_{(l,q)}(\hat
b).\tag4.20
$$

In conclusion let us emphasize that for our generalization
of the Berezin integration  the invariance of volume form (4.16) with
respect to transformation by matrix from $\s$ is essential.  It is
interest to construct the rules of  the $GL_q(n,C)$-covariant Berezin
integration.  The results of research of this point
 we intend to
publish in the following paper.

\bigskip
\centerline{\bf 5. Acknowledgements}
\bigskip
We would like to thank D.V.Volkov,  V.P.~Akulov,
 A.U.Klimyk, A.~Pashnev, V.D.~Gershun, A.V.~Mishchenko
   for critical remarks.

 V.~Rubtsov and V.~Shadura are
grateful to  Prof. A.~Niemi
for the hospitality and the exellent conditions at the Department
of Theoretical Physics of the  University of Uppsala which
made it possible to cary out this work.

We would like to thank the
members of Bogoluybov Institute for Theoretical Physics
 of the  National Academy of Science of Ukraine for the
stimulating discussions.

The work is
partially supported by grant INTAS-93-1038 (A.Bugrij,
\hfill\ \linebreak V.Shadura) and RFFI-MF0000(V.Rubtsov).
\vfill\eject

\centerline{\bf References }
\bigskip
\item{1.}
Schultz T.D., Mattis D.C., Lieb E.H. {\it Two-dimensional Ising Model
as a Soluble Problem of Many Fermions}.-- Rev. Mod. Phys., 1964, {\bf
36},  p.856-867.

McCoy B.M., Wu T.T. {\it The Two-dimensional Ising Model}\ . Harvard
University Press, Cambridge, 1973.
\item{2.} Bugrij A.I. {\it Fermionization of a
Generalized Two-Dimensional Ising \hfill\ \linebreak Model }.-- Preprint
of Institute for Theoretical Physics, ITP-86-6E, 1986, Kiev, 24p.
\item{3.}
Bugrij~A.I., Shadura~V.N. {$q$\it -Fermionization
 of the 2D Ising Model }.-- Prep-\-rint of Institute for Theoretical
 Physics, ITP-94-15R, 1994, Kiev, 38p.
 \item{ 4.} Belavin A.A., Polyakov
A.M., Zamolodchikov A.B. {\it Infinite Conformal \linebreak Symmetry in
 Two-dimensional Quantum Field Theory}.-- Nucl.   Phys. B,
1984, {\bf B241}, p.333-380.
\item{5.} Zamolodchikov A.B. {\it
Integrable Field Theory from Conformal Field.}  --
Advanced Studies in Pure Mathematics, 1989, {\bf 19}, p.641, Kinokyniya,
Tokyo.
\item{6.}
Itoyama H., Thaker H.B. {\it
Integrability and Virasoro Symmetry of the Noncritical Baxter/Ising
Model }.-- Nucl. Phys. B,
 1989, {\bf B320}, p.541.
\item{7.}
Davis B.,
Foda O., Jimbo M., Miwa T., Nakayashiki A. {\it Diagonalization of
the XXZ Hamiltonian by Vertex Operators}.-- Commun. Math. Phys., 1993,
{\bf 151}, p.89-153.
\item{8.}
Foda O., Miwa T. {\it Corner Transfer Matrices and Quantum Affine
Algeb-\-ras}.-- Int. J. Mod. Phys., 1992, {\bf A7}, sup. 1A,
p.279-302.
\item{9.} Jimbo M., Miwa T., Ohta Y. {\it Structure of
the Space of States in RSOS Models}.-- Int. J. Mod. Phys., 1993, {\bf
 A8}, p.1457-1477.
\item{10.} Foda O., Jimbo M., Miwa T., Miki K., Nakayashiki A. {\it
Vertex Operators in Solvable Lattice Models}.-- RIMS preprint
RIMS-922, Kyoto University, Kyoto, 1993, 42p.
\item{11.} Drinfeld V.G. {\it Quantum Groups}.-- Proc.Intern. Congress
of Mathematics, \linebreak Berkeley, 1986, vol.1, p.798-820.
\item{12.}
Manin Yu.I. {\it Quantum Group and Non-Commutative Geometry} . --
\linebreak Montreal University preprint CRM-1561, Montreal, 1988.
\item{13.} Faddeev L.D., Reshetikhin N.Yu., Takhtajan L.A. {\it
Quantization of Lie Qroups and Lie Algebras}.-- Algebra i Analis,
1989, {\bf 1 }, v.1, p.178-206.
\item{14.} Majid S., Marcl M.  {\it Gluing Operation for R-matrices,
Quantum Groups and Link-invariant of Hecke Type}.-- DAMTP preprint,
DAMTP/93-20, Cambridge, 1993, 36p.
\item{15.}
Lyubashenko V., Sudbery A. {\it Quantum Supergroups of $GL(n|m)$ type:
Differential Forms, Koszul Complex and Berezinians}.-- Preprint of
University of York, Heslington, 1993, 55p.
\item{16.} Wess J., Zumino
B. {\it Covariant Differential Calculus on the Quantum Hyper-\-plane}.--
Nucl. Phys. (Proc.  Suppl.), 1990, {\bf 18B}, p.302.
\item{17.} Zumino B. {\it Deformation of the Quantum Mechanical Phase
Space with Bosonic or Fermionic Coordinates}.-- Mod. Phys. Lett. A,
1991, {\bf A 13}, p.1225.
\item{18.} Woronowicz S.L. {\it Differential Calculus on Compact
Matrix Pseudogroups} .-- Commun. Math. Phys., 1989, {\bf 122},
p.125-170.

\item{19.}
Gurevich D., Radul A., Rubtsov V. {\it Noncommutative Differential
Geometry Connected with Yang-Baxter Equation}.-- Zap. Nauch. Sem.
LOMI, vol.199, 1992, p.51-70.
\item{20}
Maltsiniotis G..{\it Le Langage des Espases et des Groupes Quantiques}.
--Commun.Math.Phys. 1993,{\bf 151},p.275-302.
\item{21.}
Berezin F.A. {\it Introduction to Algebra and Analysis with
Anticommutative Variables }. Published by Moscow University Press,
 Moscow, 1983, 205p.
 \item{22.}
 Baulieu L.,
Floratos E.G. {\it Path Integral on the Quantum Plane}.-- Phys. Lett.
B, 1991, {\bf B258}, p.171.
\item{23.}
Chachian M., Demichev A.P.
{\it q-Deformed Path Integral}.-- Preprint of Hel-\-sinki University,
 HU - SEFT R 1993-10, Helsinki, 1993, 14p.
 \item{24.}
 Mathai
V., Quillen D. {\it Superconnections, Thom Classes and Equivariant
Differential Forms}.-- Topology, 1986, {\bf 25}, p.85-110.

\end